\RequirePackage{fix-cm}
\documentclass[envcountsect, envcountsame]{svjour3}                     % onecolumn (standard format)

\smartqed  % flush right qed marks, e.g. at end of proof

\usepackage{cite}

\usepackage[x11names]{xcolor}
\usepackage[
  colorlinks,
]{hyperref}

\newcommand\rurl[1]{%
  \href{http://#1}{\nolinkurl{#1}}%
}

\AtBeginDocument{\hypersetup{citecolor=magenta, urlcolor=magenta, linkcolor=blue}}

\usepackage{lmodern}
\usepackage[T5]{fontenc}

\usepackage{pb-diagram}

\usepackage{amsmath}
\usepackage{bm}
\usepackage{amssymb}
\usepackage{mathrsfs}
\usepackage{latexsym}
\usepackage{exscale}
\usepackage{eucal}

\usepackage{listings}

\def\DD{D\kern-.7em\raise0.4ex\hbox{\char '55}\kern.33em}

\usepackage[top=2.5cm, bottom=2.5cm, left=2cm, right=2cm] {geometry}
\def\subclassname{{\bfseries Mathematics Subject Classification
(2010)}\enspace}
\def\subclass#1{\par\addvspace\medskipamount{\rightskip=0pt plus1cm
\def\and{\ifhmode\unskip\nobreak\fi\ $\cdot$
}\noindent\subclassname\ignorespaces#1\par}}

\spnewtheorem{thm}{Theorem}[section]{\bfseries}{\itshape}
\spnewtheorem{prop}[thm]{Proposition}{\bfseries}{\itshape}
\spnewtheorem{lem}[thm]{Lemma}{\bfseries}{\itshape}
\spnewtheorem{gt}[thm]{Conjecture}{\bfseries}{\itshape}
\spnewtheorem{corl}[thm]{Corollary}{\bfseries}{\itshape}

\spnewtheorem{defn}[thm]{Definition}{\bfseries}{\upshape}
\spnewtheorem{rem}[thm]{Remark}{\bfseries}{\upshape}
\spnewtheorem{nx}[thm]{Remark}{\bfseries}{\upshape}
\spnewtheorem{nota}[thm]{Notation}{\bfseries}{\upshape}

\spnewtheorem{therm}{Theorem}[section]{\bfseries}{\itshape}
\spnewtheorem{propo}[therm]{Proposition}{\bfseries}{\itshape}
\spnewtheorem{lema}[therm]{Lemma}{\bfseries}{\itshape}
\spnewtheorem{conj}[therm]{Conjecture}{\bfseries}{\itshape}
\spnewtheorem{conje}[therm]{Conjecture}{\bfseries}{\itshape}
\spnewtheorem{corls}[therm]{Corollary}{\bfseries}{\itshape}

\spnewtheorem{defns}[therm]{Definition}{\bfseries}{\upshape}
\spnewtheorem{rems}[therm]{Remark}{\bfseries}{\upshape}
\spnewtheorem{ex}[thm]{Example}{\bfseries}{\upshape}
\spnewtheorem{notas}[therm]{Notation}{\bfseries}{\upshape}

\spnewtheorem{dn}{Definition}[subsection]{\bfseries}{\upshape}
\spnewtheorem{kh}[dn]{Notation}{\bfseries}{\upshape}
\spnewtheorem{nte}[thm]{Note}{\bfseries}{\upshape}

\spnewtheorem{bd}[dn]{Lemma}{\bfseries}{\itshape}
\spnewtheorem{theo}[dn]{Theorem}{\bfseries}{\itshape}
\spnewtheorem{mde}[dn]{Proposition}{\bfseries}{\itshape}

\spnewtheorem*{acknow}{Acknowledgment}{\bfseries}{\upshape}

\spnewtheorem*{xproof}{}{\itshape}{\rmfamily}% the label is assigned later

\renewenvironment{proof}[1][\proofname]
 {\renewcommand\xproofname{#1}\xproof}
 {\endxproof}

\def\DD{D\kern-.7em\raise0.4ex\hbox{\char '55}\kern.33em}

\begin{document}

%\title{On generators of the unstable $\mathscr A$-module $H^{*}((K(\mathbb F_2, 1))^{\times 5})$\\ in a generic degree and applications\\[2mm]{\fontsize{9pt}{9pt}\selectfont \textit{(Dedicated to Prof. William Singer on the occasion of his 76th birthday)}}}

%\titlerunning{$\mathscr A$-generators of $H^{*}((K(\mathbb F_2, 1))^{\times 5})$ in a generic degree and applications}        % if too long for running head

\title{On modules over the mod 2 Steenrod algebra and hit problems}
%{\fontsize{9pt}{9pt}\selectfont \textit{(Dedicated to Prof. Masaki Kameko, Prof. Ali S. Janfada, and Prof. William Singer)}}}

\titlerunning{On modules over the mod 2 Steenrod algebra and hit problems}        % if too long for running head

\author{\DD\d{\u a}ng V\~o Ph\'uc
         %etc.
}

\authorrunning{\DD\d.V.Ph\'uc} % if too long for running head

\institute{\DD\d{\u a}ng V\~o Ph\'uc \at
        Faculty of Education Studies\\
 University of Khanh Hoa, Nha Trang, Khanh Hoa, Vietnam \\
              \email{dangphuc150488@gmail.com}  \\
               %This work is supported by National Foundation for Science and Technology Development\\ (NAFOSTED) of Viet Nam [grant number 101.04-2017.05]
}
         %  \\
%             \emph{Present address:} of F. Author  %  if needed

%\date{Received: date / Accepted: date}
% The correct dates will be entered by the editor

\maketitle

\begin{abstract}

Let us consider the prime field of two elements, $\mathbb F_2\equiv \mathbb Z_2.$ It is well-known that the classical "hit problem" for a module over the mod 2 Steenrod algebra $\mathscr A$ is an interesting and important open problem of Algebraic topology, which asks a minimal set of generators for the polynomial algebra $\mathcal P_m:=\mathbb F_2[x_1, x_2, \ldots, x_m]$, regarded as a connected unstable $\mathscr A$-module on $m$ variables $x_1, \ldots, x_m,$ each of degree 1. The algebra $\mathcal P_m$ is the $\mathbb F_2$-cohomology of the product of $m$ copies of the Eilenberg-MacLan complex $K(\mathbb F_2, 1).$ Although the hit problem has been thoroughly studied for more than 3 decades, solving it remains a mystery for $m\geq 5.$  It is our intent in this work is of studying the hit problem of five variables. More precisely, we develop our previous work [Commun. Korean Math. Soc. 35 (2020), 371-399] on the hit problem for $\mathscr A$-module $\mathcal P_5$ in a degree of the generic form $n_t:=5(2^t-1) + 18.2^t,$ for any non-negative integer $t.$ An efficient approach to solve this problem had been presented. 
%As a consequence, the calculations confirmed Sum's conjecture \cite{N.S2} for the relationship between the minimal sets of $\mathscr A$-generators of the polynomial algebras $\mathcal P_{m-1}$ and $\mathcal P_{m}$ in the case $m=5$ and degree $n_t.$ 
Two applications of this study are to determine the dimension of $\mathcal P_6$ in the generic degree $5(2^{t+4}-1) + n_1.2^{t+4}$ for all $t > 0$ and to describe the modular representations of the general linear group of rank 5 over $\mathbb F_2.$ As a corollary, the cohomological "transfer", defined by William Singer [Math. Z. 202 (1989), 493-523], is an isomorphism in bidegree $(5, 5+n_0).$ Singer's transfer is one of the relatively efficient tools to approach the structure of mod-2 cohomology of the Steenrod algebra.
%Moreover, we provide an algorithm in MAGMA for verifying the results and studying the hit problem in general. 
%The results of this paper are the special cases but they provide useful information in formulating general solutions. 

 \subclass{ 13A50\and 55Q45\and 55S10\and 55S05\and 55T15\and 55R12}
\end{abstract}

\section{Introduction}

Let $\mathcal{O}^S(i, \mathbb{F}_2, \mathbb{F}_2)$ denote the set of all stable cohomology operations of degree $i,$ with coefficient in the prime field $\mathbb F_2.$ Then, the $\mathbb F_2$-algebra $\mathscr A:=\bigoplus_{i\geq 0}\mathcal{O}^S(i, \mathbb{F}_2, \mathbb{F}_2)$ is called \textit{the mod 2 Steenrod algebra}. In other words, the algebra $\mathscr A$ is the algebra of stable operations on the mod 2 cohomology.  In \cite{J.M}, Milnor observed that this algebra is also a graded connected cocommutative Hopf algebra over $\mathbb F_2.$ %Therefore, its dual $\mathscr A^{*}$ is a commutative algebra and in the same paper, Milnor showed that $\mathscr A^{*}$ is a polynomial algebra. 
In some cases, the resulting $\mathscr A$-module structure on $H^{*}(X, \mathbb F_2)$ provides additional information about CW-complexes $X;$ for instance (see section three for a detailed proof), the CW-complexes $\mathbb{C}P^4/\mathbb{C}P^2$ and $\mathbb{S}^6\vee \mathbb{S}^8$ have cohomology rings that agree as a graded commutative $\mathbb F_2$-algebras, but are different as a module over $\mathscr A.$ Afterwards, the Steenrod algebra is widely studied by mathematicians whose interests range from algebraic topology and homotopy theory to manifold theory, combinatorics, representation theory, and more. It is well-known that the $\mathbb F_2$-cohomology of the Eilenberg-MacLan complex $K(\mathbb F_2, 1)$ is isomorphic to $\mathbb F_2[x],$ the polynomial ring of degree $1$ in one variable. Hence, based upon the K\"{u}nneth formula for cohomology, we have an isomorphism of $\mathbb F_2$-algebras
$$\mathcal P_m:=H^*((K(\mathbb F_2, 1))^{\times m}, \mathbb F_2) \cong \underset{\text{ $m$ times}}{\underbrace{ \mathbb{F}_2[x_1]\otimes_{\mathbb{F}_2}  \mathbb{F}_2[x_2]\otimes_{\mathbb{F}_2}\cdots\otimes_{\mathbb{F}_2} \mathbb{F}_2[x_m]}}\cong \mathbb{F}_2[x_1, \ldots, x_m],$$ where $x_i\in H^1((K(\mathbb F_2, 1))^{\times m}, \mathbb F_2)$ for every $i.$ Since $\mathcal P_m$ is the cohomology of a CW-complex, it is equipped with a structure of unstable module over $\mathscr A.$ It has been known (see also \cite{S.E}) that $\mathscr A$ is spanned by the Steenrod squares $Sq^{i}$ of degree $i$ for $i\geq 0$ and that the action of $\mathscr{A}$ on $\mathcal P_m$ is depicted as follows:
$$ \begin{array}{ll}
Sq^i(x_t) &= \left\{\begin{array}{lll}
x_t &\mbox{if}& i = 0,\\
x_t^2 &\mbox{if}& i = 1, \ \ \mbox{(\textit{the instability condition})},\\
0 &\mbox{if}& i > 1,
\end{array}\right.\\
Sq^i(FG) &= \sum_{0\leq \alpha\leq i}Sq^{\alpha}(F)Sq^{i-\alpha}(G),\ \mbox{for all $F,\ G\in \mathcal P_m$}\ \ (\mbox{\textit{the Cartan formula})}.
\end{array}$$
It is to be noted that since $Sq^{\deg(F)}(F) = F^{2}$ for any $F\in \mathcal P_m$, the polynomial ring $\mathcal P_m$ is also an unstable $\mathscr A$-algebra. Letting $GL_m:=GL(m, \mathbb F_2)$ for the general linear group of degree $m$ over $\mathbb F_2.$ This $GL_m$ when $m\geq 2,$ which can be generated by two elements (see Waterhouse \cite{Waterhouse}), acts on $\mathcal P_m$ by matrix substitution. So, in addition to $\mathscr A$-module structure, $\mathcal P_m$ is also a (right) $\mathbb F_2GL_m$-module. The classical "hit problem" for the algebra $\mathscr A$, which is concerned with seeking a minimal set of $\mathscr A$-generators for $\mathcal P_m$, has been initiated in a variety of contexts by Peterson \cite{F.P}, Priddy \cite{S.P}, Singer \cite{W.S2}, and Wood \cite{R.W}. Structure of modules over $\mathscr A$ and hit problems are currently one of the central subjects in Algebraic topology and has a great deal of intensively studied by many authors like Brunetti and collaborators \cite{Brunetti1, Brunetti2}, Crabb-Hubbuck \cite{C.H}, Inoue \cite{M.I1, M.I2}, Janfada-Wood \cite{J.W, J.W2}, Janfada \cite{Janfada1, Janfada2}, Kameko \cite{M.K}, Mothebe-Uys \cite{M.M}, Mothebe \cite{M.M2}, Pengelley-William \cite{P.W}, the present author and N. Sum \cite{P.S1, P.S2, D.P1, D.P2, D.P3, D.P6, D.P7, N.S1, N.S2, N.S3}, Walker-Wood \cite{W.W, W.W2}, etc. As it is known, when $\mathbb F_2$ is an $\mathscr A$-module concentrated in degree 0, solving the hit problem is to determine an $\mathbb F_2$-basis for the space of indecomposables, or "unhit" elements, $Q^{\otimes m}:=  \mathbb{F}_2\otimes_{\mathscr{A}} \mathcal P_m = \mathcal P_m/ \overline{\mathscr A}\mathcal P_m$ where $\overline{\mathscr A}$ is the positive degree part of $\mathscr A$. It is well-known that the action of $GL_m$ and the action of $\mathscr A$ on $\mathcal P_m$ commute. So, there is an induced action of  $GL_m$ on $Q^{\otimes m}.$ The structure of $Q^{\otimes m}$ has been treated for $m\leq 4$ by Peterson \cite{F.P}, Kameko \cite{M.K} and Sum \cite{N.S1}. The general case is an interesting open problem. Most notably, the study of this space plays a vital role in describing the $E^{2}$-term of the Adams spectral sequence (Adams SS), ${\rm Ext}_{\mathscr{A}}^{m,m+*}(\mathbb{F}_2, \mathbb{F}_2)$ via the $m$-th Singer cohomological "transfer" \cite{W.S1}. This transfer is a linear map $$Tr^{\mathscr A}_m: (\mathbb F_2 \otimes_{GL_m}P_{\mathscr A}((\mathcal P_m)^{*}))_n\to {\rm Ext}_{\mathscr {A}}^{m,m+n}(\mathbb F_2,\mathbb F_2) = H^{m, m+n}(\mathscr A, \mathbb F_2),$$ from the subspace of all $\overline{\mathscr A}$-annihilated elements to the $E^2$-term of the Adams SS. Here $(\mathcal P_m)^{*} = H_*((K(\mathbb F_2, 1))^{\times m}, \mathbb F_2) $ and $\mathbb F_2 \otimes_{GL_m}P_{\mathscr A}((\mathcal P_m)^{*})$ are the dual of $\mathcal P_m$ and $(Q^{\otimes m})^{GL_m},$ respectively, where $(Q^{\otimes m})^{GL_m}$ denotes the space of $GL_m$-invariants.
%The algebra $(\mathcal P_m)^{*}$ has a right $\mathscr A$-module structure. The right action of $\mathscr A$ on this algebra is given by $$(a^{(j)}_t)Sq^{k} = \binom{j-k}{k}a^{(j-k)}_t= Sq_*^{k}(a^{(j)}_t)$$ and Cartan's formula, where $Sq_*^{k}$ denotes the dual Steenrod squares. 
A natural question arises: Why do we need to calculate the Adams $E^2$-term? The answer is that it is involved in determining the stable homotopy groups of spheres. These groups are pretty fundamental and interesting. Nevertheless, they are also not fully-understood subjects yet. Therefore, the clarification of these problems is an important task of Algebraic topology. It has been shown (see \cite{J.B}, \cite{W.S1}) that the algebraic transfer is highly nontrivial, more precisely, that $Tr^{\mathscr A}_m$ is an isomorphism for $0 < m < 4$ and that the "total" transfer $\bigoplus_{m\geq 0}Tr^{\mathscr A}_m: \bigoplus_{m\geq 0}(\mathbb F_2 \otimes_{GL_m}P_{\mathscr A}((\mathcal P_m)^{*}))_n\to \bigoplus_{m\geq 0}{\rm Ext}_{\mathscr {A}}^{m,m+n}(\mathbb F_2,\mathbb F_2)$ is a homomorphism of bigraded algebras with respect to the product by concatenation in the domain and the usual Yoneda product for the Ext group. Minami's works \cite{N.M1, N.M2} have shown the usefulness of the Singer transfer and the hit problem for surveying the \textit{Kervaire invariant one problem}. This problem, which is a long standing open topic in Algebraic topology, asks when there are framed manifolds with Kervaire invariant one. (Note that a \textit{framing} on a closed smooth manifold $M^{n}$ is a trivialization of the normal bundle $\nu(M, i)$ of some smooth embedding  $i: M\hookrightarrow \mathbb R^{n+*}.$ Here $\nu(M, i)$ is defined to be a quotient of the pullback of the tangent bundle of $\mathbb R^{n+*}$ by the sub-bundle given by the tangent bundle of $M.$ So, $\nu(M, i)$ is an $*$-dimensional real vector bundle over $M^{n}.$ For more details, we refer the reader to \cite{Snaith}.) Framed manifolds of Kervaire invariant one have been constructed in dimension $2^{k}-2$ for $2\leq k\leq 6.$ In 2016, by using mod 8 equivariant homotopy theory, Hill, Hopkins, and Ravenel claimed in their surprising work \cite{Hill} that the Kervaire invariant is 0 in dimension $2^{k}-2$ for $k\geq 8.$ Up to present, it remains undetermined for $k = 7$ (or dimension $126$) and this has the status of a hypothesis by Snaith \cite{Snaith}. 

Return to Singer's transfer, in higher homological degrees, the works \cite{B.H.H}, \cite{Ha}, \cite{Hung}, \cite{T.N2}, and \cite{H.Q} determined completely the image of $Tr_4^{\mathscr A}$. The authors show that the image of the fourth transfer contains every element in the four families $\{d_t|\, t\geq 0\},\ \{e_t|\, t\geq 0\}$, $\{f_t|\, t\geq 0\}$, and $\{p_t|\, t\geq 0\}$, whereas it does not contain any element in the three families $\{g_{t+1}|\, t\geq 0\}$, $\{D_3(t)|\, t\geq 0\}$, and $\{p'_t|\, t\geq 0\}$. More explicitly, the result on $\{g_{t+1}|\, t\geq 0\}$ is due to \cite{B.H.H}; that on $\{D_3(t)|\, t\geq 0\}$, and $\{p'_t|\, t\geq 0\}$ is due to \cite{Hung}; the conclusion on $\{d_t|\, t\geq 0\},$ and $\{e_t|\, t\geq 0\}$ is proved by \cite{Ha}, while that on $\{f_t|\, t\geq 0\}$ is proved by \cite{T.N2}. Remarkably, the results by \cite{B.H.H} and \cite{Hung} gave a negative answer to Minami's hypothesis \cite{N.M2} predicting that the localization of $Tr_m^{\mathscr A}$ given by inverting the squaring operation $Sq^{0}$ (see section two) is an isomorphism. In \cite{Hung}, H\uhorn ng indicated that $Tr_4^{\mathscr A}$ is not an isomorphism in infinitely many degrees. In particular, from preliminary calculations in \cite{W.S1}, Singer proposed the following.

\begin{gt}\label{gtSin}
The transfer homomorphism is a monomorphism in  every rank $m > 0.$
\end{gt}

One has seen above that $Tr_m^{\mathscr A}$ is an isomorphism for $m<4,$ and so the conjecture holds in these ranks $m.$ Our recent work \cite{D.P14} has shown that it is also true for $m = 4,$ but the answer to the general case remains a mystery, even in the case of $m=5$ with the help of a computer algebra. It is known, in ranks $\leq 4$, the calculations of Singer \cite{W.S1}, H\`a \cite{Ha}, and Nam \cite{T.N2} tell us that the non-zero elements $h_t\in {\rm Ext}_{\mathscr {A}}^{1, 2^{t}}(\mathbb F_2,\mathbb F_2),\, e_t\in {\rm Ext}_{\mathscr {A}}^{4,2^{t+4} + 2^{t+2} + 2^{t}}(\mathbb F_2,\mathbb F_2),\, f_t\in {\rm Ext}_{\mathscr {A}}^{4,2^{t+4} + 2^{t+2} + 2^{t+1}}(\mathbb F_2,\mathbb F_2),$ for all $t\geq 0,$ are detected by the cohomological transfer. In rank 5, based on invariant theory, Singer \cite{W.S1} gives an explicit element in ${\rm Ext}_{\mathscr {A}}^{5,5+9}(\mathbb F_2,\mathbb F_2),$ namely $Ph_1$, that is not detected by $Tr_5^{\mathscr A}.$ In general, direct calculating the value of $Tr_m^{\mathscr A}$ on any non-zero element is difficult. Moreover, there is no general rule for that, and so, each computation is important on its own. By this and the above results, in the present text, we would like to investigate the family $\{h_tf_t = h_{t+1}e_t\in {\rm Ext}_{\mathscr {A}}^{5,5+(23.2^{t}-5)}(\mathbb F_2,\mathbb F_2)|\, t\geq 0\},$ and Singer's conjecture for $m=5$ in degree $ 5(2^{t}-1) + 18.2^{t}=23.2^{t}-5$ with $t = 0.$ To do this, we use a basis of the indecomposables $Q^{\otimes 5}$ in degree $18 = 5(2^{0}-1) + 18.2^{0}$, which is given by our previous work \cite{D.P2} (see Proposition \ref{md18} below). In addition, the main goal of this work is to also compute explicitly the dimension of $Q^{\otimes 5}$ in degree $5(2^{t}-1) + 18.2^{t}$ for the cases $t\geq 1.$ Then, Singer's conjecture for $m=5$ and these degrees will be discussed at the end of section two. We hope that our results would be necessary to formulate general solutions.  %This approach is different from the ones of Singer \cite{W.S1}.

\section{Statement of results}

%The hit problem was solved explicitly by Peterson \cite{F.P}, Kameko \cite{M.K} and Sum \cite{N.S1} for the number of variables $\leq 4.$ 
{\bf Some notes.} %We develop our previous result in \cite{D.P3} on the structure of $Q^{\otimes 5}$ in generic degree $5(2^t-1) + 18.2^t$ for any non-negative integer $t.$ As an application, we examine Conjecture \ref{gtSin} for $Tr_5^{\mathscr A}$ in the above degree for $t = 0.$ The result when $t\geq 1$ was used to determine the dimension of $Q^{\otimes 6}$ in the generic degree $5(2^{t+4}-1) + 41.2^{t+4}.$ 
Throughout this paper, let us write 
$$ \begin{array}{ll}
\medskip
(\mathcal P_m)_n&:= \langle \{f\in \mathcal P_m|\, \mbox{$f$ is a homogeneous polynomial of degree $n$}\} \rangle,\\
Q^{\otimes m}_n&:= \langle \{[f]\in Q^{\otimes m}|\, f\in (\mathcal P_m)_n\}\rangle,
\end{array}$$ which are $\mathbb F_2GL_m$-submodules of $\mathcal P_m$ and $Q^{\otimes m},$ respectively. So $\mathcal P_m = \bigoplus_{n\geq 0}(\mathcal P_m)_n$ and $Q^{\otimes m} = \bigoplus_{n\geq 0}Q^{\otimes m}_n.$ Recall that to solve the hit problem of three variables, Kameko \cite{M.K} constructed a $\mathbb F_2GL_m$-modules epimorphism:
$$ \begin{array}{ll}
(\widetilde {Sq^0_*})_{(m, m+2n)}: Q^{\otimes m}_{m+2n}  &\longrightarrow Q^{\otimes m}_n\\
\ \ \ \ \ \ \ \ \ \ \ \ \ \ \ \ \mbox{[}\prod_{1\leq j\leq m}x_j^{a_j}\mbox{]}&\longmapsto \left\{\begin{array}{ll}
\mbox{[}\prod_{1\leq j\leq m}x_j^{\frac{a_j-1}{2}}\mbox{]}& \text{if $a_j$ odd, $j = 1, 2, \ldots, m$},\\
0& \text{otherwise},
\end{array}\right.
\end{array}$$ 
which induces the homomorphism $\widetilde {Sq^0_*}: (Q^{\otimes m}_{m+2n})^{GL_m}\to (Q^{\otimes m}_{n})^{GL_m}.$   Since $\mathscr A$ is a cocommutative Hopf algebra, there exists the squaring operations $Sq^i: {\rm Ext}_{\mathscr A}^{m, m+n}(\mathbb F_2, \mathbb F_2) \to {\rm Ext}_{\mathscr A}^{m+i, 2m+2n}(\mathbb F_2, \mathbb F_2),$ which share most of the properties with $Sq^{i}$ on the cohomology of spaces (see \cite{May}), but the classical $Sq^{0}$ is not the identity in general.
%It should be note that  $Sq^{0}$ is a special case of the Steenrod operations $Sq^{i}: {\rm Ext}_{\mathscr A}^{s, t}(\mathbb U, \mathbb V)\to {\rm Ext}_{\mathscr A}^{s+i, 2t}(\mathbb U, \mathbb V),$ where $\mathbb U$ is an $\mathscr A$-coalgebra and $\mathbb V$ is an $\mathscr A$-algebra. 
Remarkably, this $Sq^{0}$ commutes with the dual of $\widetilde {Sq^0_*}$ through the Singer transfer (see \cite{J.B}, \cite{N.M2}). The reader who is familiar with Kameko's $(\widetilde {Sq^0_*})_{(m, m+2n)}$ will probably agree that this map is very useful in solving the hit problem. Indeed, Kameko \cite{M.K} showed that if  $m = \xi(n) = {\rm min}\{\gamma\in\mathbb N:\ n = \sum_{1\leq i\leq \gamma}(2^{d_i}-1),\, d_i > 0, \forall i,\, 1\leq i\leq \gamma\},$ then $(\widetilde {Sq^0_*})_{(m, m+2n)}$ is an isomorphism of $\mathbb F_2GL_m$-modules. This statement and Wood's work \cite{R.W} together are sufficient to determine $Q^{\otimes m}_n$ in each degree $n$ of the special "generic" form $n = r(2^t-1) + d.2^t,$ whenever $0< \xi(d) < r < m,$ and $t\geq 0$ (see also \cite{D.P6}). 

As we mentioned at the beginning, the hit problem was completely solved for $m\leq 4.$ Very little information is known for $m = 5$ and degrees $n$ given above. At least, it is surveyed by the present writer \cite{D.P6} for $(r, d, t)\in \{(5, 18, 0), (5, 8, t)\}.$ We now extend for the case $(r, d, t) = (5, 18, t),$ in which $t$ an arbitrary non-negative integer. We start with a useful remark.

\begin{rem}\label{nx1}
It can be easily seen that $5(2^t-1) + 18.2^t = 2^{t+4} + 2^{t+2} + 2^{t+1} + 2^{t-1} + 2^{t-1} - 5,$ and so $\xi(5(2^t-1) + 18.2^t) = 5$ for any $t > 1.$ This implies that the iterated Kameko map $$((\widetilde {Sq^0_*})_{(5, 5(2^{t} - 1) + 18.2^{t})})^{t-1}:  Q^{\otimes 5}_{5(2^{t} - 1) + 18.2^{t}} \to Q^{\otimes 5}_{5(2^{1} - 1) + 18.2^{1}}$$ is an isomorphism, for all $t\geq 1,$ and therefore, it is enough to determine $Q^{\otimes 5}_{5(2^{t} - 1) + 18.2^{t}}$ for $t\in \{0, 1\}$. The case $t = 0$ has explicitly been computed by us in \cite{D.P3}. When $t = 1,$ because Kameko's homomorphism $$(\widetilde {Sq^0_*})_{(5, 5(2^{1} - 1) + 18.2^{1})}: Q^{\otimes 5}_{5(2^{1} - 1) + 18.2^{1}} \to Q^{\otimes 5}_{5(2^{0} - 1) + 18.2^{0}}$$ is an epimorphism, we have an isomorphism
$$ Q^{\otimes 5}_{5(2^{1} - 1) + 18.2^{1}}  \cong {\rm Ker}((\widetilde {Sq^0_*})_{(5, 5(2^{1} - 1) + 18.2^{1})})\bigoplus Q^{\otimes 5}_{5(2^{0} - 1) + 18.2^{0}}.$$ 
The space $Q^{\otimes 5}_{5(2^{0} - 1) + 18.2^{0}}$ is known by our previous work \cite{D.P3}. Thus, we need compute the kernel of $(\widetilde {Sq^0_*})_{(5, 5(2^{1} - 1) + 18.2^{1})}.$  For this, our approach can be summarized as follows:

\begin{itemize}

\item[(i)] A mononial in $\mathcal P_5$ is assigned a \textit{weight vector} $\omega$ of degree $5(2^{1} - 1) + 18.2^{1}$, which stems from the binary expansion of the exponents of the monomial. The space of indecomposable elements ${\rm Ker}((\widetilde {Sq^0_*})_{(5, 5(2^{1} - 1) + 18.2^{1})})$ is then decomposed into a direct sum of $(Q_{5(2^{1} - 1) + 18.2^{1}}^{\otimes 5})^{0}$ and the subspaces $(Q^{\otimes 5})^{\omega^{>0}}$ indexed by the weight vectors $\omega.$ Here $[F]_{\omega} = [G]_{\omega}$ in $(Q_{5(2^{1} - 1) + 18.2^{1}}^{\otimes 5})^{\omega}$ if the polynomial $F-G$ is hit, modulo a sum of monomials of weight vectors less than $\omega.$ Basing the previous results by Peterson \cite{F.P}, Kameko \cite{M.K}, Sum \cite{N.S1}, and by us \cite{D.P6}, one can easily determine $(Q_{5(2^{1} - 1) + 18.2^{1}}^{\otimes 5})^{0}.$

\item[(ii)] The monomials in a given degree are lexicographically ordered first by weight vectors and then by exponent vectors. This leads to the concept of \textit{admissible monomial}; more explicitly, a monomial is admissible if, modulo hit elements, it is not equal to a sum of monomials of smaller orders. The space $(Q_{5(2^{1} - 1) + 18.2^{1}}^{\otimes 5})^{\omega^{>0}}$ above is easily seen to be isomorphic to the space generated by admissible monomials of the weight vector $\omega.$

\item[(iii)] In a given (small) degree, we first list all possible weight vectors of an admissible monomial. This is done by first using a criterion of Singer \cite{W.S1} on the hit monomials, and then combining with the results by Kameko \cite{M.K} and Sum \cite{N.S1} (see Theorems \textcolor[rgb]{0.0,0.0,1.0}{4.2}, and \textcolor[rgb]{0.0,0.0,1.0}{4.3} in section four)
%(see Theorems \ref{dlK} and \ref{dlS} in section four) 
of the form "$XZ^{2^{r}}$ (or $ZY^{2^{t}}$) admissible implying $Z$ admissible, under some mild conditions".

\item[(iv)] In a given weight vector, we claim the (strict) inadmissibility of some explicit monomials. The proof is given for a typical monomial in each case by explicit computations. Finally, a direct calculation using %Theorems \ref{dlMM}, \ref{dlbs}, 
Theorems  \textcolor[rgb]{0.0,0.0,1.0}{3.2}, \textcolor[rgb]{0.0,0.0,1.0}{3.3}, and some homomorphisms in section three, we obtain a basis of $(Q_{5(2^{1} - 1) + 18.2^{1}}^{\otimes 5})^{\omega^{>0}}.$ This approach is much less computational and it can be applied for all certain degrees and all variables $m.$ Moreover, the MAGMA computer algebra \cite{Magma} has been used for verifying the results. 
%We are thankful to Prof. Bob Bruner for sending this program. 
%experimentation leading up to some of the results below. 
%In addition, to verify the results, we have used an algorithm on the MAGMA computer algebra system \cite{Magma} written by Robert R. Bruner. We are thankful to Professor Bruner for sending this program. 
\end{itemize}
\end{rem}

Before going into detail and  proceeding to the main results, let us provide some basic concepts. Of course, we assume that the reader is not familiar with the
basics of hit problems.
%in order to state our main result precisely
\medskip

{\bf Weight vector and exponent vector}. Let $\omega = (\omega_1, \omega_2, \ldots, \omega_t,\ldots)$ be a sequence of non-negative integers. We say that $\omega$ is a \textit{weight vector}, if $\omega_t  = 0$ for $t\gg 0.$ Then, we also define $\deg(\omega) = \sum_{t\geq 1}2^{t-1}\omega_t.$ Let $X = x_1^{u_1}x_2^{u_2}\ldots x_m^{u_m}$ be a mononial in $\mathcal P_m,$ define two sequences associated with $X$ by $$\omega(X) :=(\omega_1(X), \omega_2(X), \ldots, \omega_t(X), \ldots),\ \ u(X):=(u_1, u_2, \ldots, u_m),$$ where $\omega_t(X)=\sum_{1\leq j\leq m}\alpha_{t-1}(u_j)$ in which $\alpha_t(n)$ denotes the $t$-th coefficients in dyadic expansion of a positive integer $n.$ They are called the {\it weight vector} and the \textit{exponent vector} of $X,$ respectively. We use the convention that the sets of all the weight vectors and the exponent vectors are given the left lexicographical order.   

\medskip

{\bf Linear order on \boldmath{$\mathcal P_m$}}. Assume that $X = x_1^{u_1}x_2^{u_2}\ldots x_m^{u_m}$ and $Y = x_1^{v_1}x_2^{v_2}\ldots x_m^{v_m}$ are the monomials of the same degree in $\mathcal P_m.$ We say that $X  < Y$ if and only if one of the following holds:
\begin{enumerate}
\item[(i)] $\omega(X) < \omega(Y);$
\item[(ii)] $\omega(X) = \omega(Y)$ and $u(X) < v(Y).$
\end{enumerate}

\medskip

{\bf Equivalence relations on \boldmath{$\mathcal P_m$}}. For a weight vector $\omega,$ we denote two subspaces associated with $\omega$ by
$$ \begin{array}{ll}
 \mathcal P^{\leq \omega}_m = \langle\{ X\in\mathcal P_m|\, \deg(X) = \deg(\omega),\  \omega(X)\leq \omega\}\rangle,\\
 \mathcal P^{<\omega}_m = \langle \{ X\in\mathcal P_m|\, \deg(X) = \deg(\omega),\  \omega(X) < \omega\}\rangle.
\end{array}$$ 
Let $F$ and $G$ be the homogeneous polynomials in $\mathcal P_m$ such that $\deg(F) = \deg(G).$ We say that 
\begin{enumerate}
\item [(i)]$F \equiv G $ if and only if $(F - G)\in \overline{\mathscr {A}}\mathcal P_m=\sum_{i\geq 0}{\rm Im}(Sq^{2^{i}}).$ Specifically,  if $F \equiv 0,$ then $F$ is hit (or $\mathscr A$-decomposable), i.e., $F$ can be written in the form $\sum_{i\geq 0}Sq^{2^{i}}(F_i)$ for some $F_i\in \mathcal P_m$;
\item[(ii)] $F \equiv_{\omega} G$ if and only if $F, \, G\in \mathcal P^{\leq \omega}_m$ and $(F -G)\in ((\overline{\mathscr {A}}\mathcal P_m \cap \mathcal P_m^{\leq \omega})+ \mathcal P_m^{<\omega}).$
\end{enumerate}
It is not difficult to show that the binary relations "$\equiv$" and "$\equiv_{\omega}$" are equivalence ones. So, one defines the quotient space 
$$(Q^{\otimes m})^{\omega} = \mathcal P_m^{\leq \omega}/ ((\overline{\mathscr {A}}\mathcal P_m \cap \mathcal P_m^{\leq \omega})+ \mathcal P_m^{<\omega}).$$
Moreover, due to Sum \cite{N.S3}, $(Q^{\otimes m})^{\omega}$ is also an $\mathbb F_2GL_m$-module. 

\medskip

{\bf Admissible monomial and inadmissible monomial}. A monomial $X\in\mathcal P_m$ is said to be {\it inadmissible} if there exist monomials $Y_1, Y_2,\ldots, Y_k$ such that $Y_j < X$ for $1\leq j\leq k$ and $X \equiv  \sum_{1\leq j\leq k}Y_j.$ Then, $X$ is said to be {\it admissible} if it is not inadmissible.

Thus, with the above definitions in hand, it is straightforward to see that the set of all the admissible monomials of degree $n$ in $\mathcal P_m$ is \textit{a minimal set of $\mathscr {A}$-generators for $\mathcal P_m$ in degree $n.$} So, $Q_n^{\otimes m}$ is a $\mathbb F_2$-vector space with a basis consisting of all the classes represent by the admissible monomials of degree $n$ in $\mathcal P_m.$ Further, as stated in \cite{D.P2}, the dimension of $Q_n^{\otimes m}$ can be represented as the sum of the dimensions $(Q^{\otimes m})^{\omega}$ such that $\deg(\omega) = n.$ For later convenience, we need to set some notation. Let $\mathcal P^{0}_m$ and $\mathcal P^{>0}_m$ denote the $\mathscr A$-submodules of $\mathcal P_m$ spanned all the monomials $\prod_{1\leq j\leq m}x_j^{t_j}$ such that $\prod_{1\leq j\leq m}t_j = 0,$ and $\prod_{1\leq j\leq m}t_j > 0,$ respectively. Let us write $(Q^{\otimes m})^0:= \mathbb F_2\otimes_{\mathscr A} \mathcal P^{0}_m,\ \mbox{and}\ (Q^{\otimes m})^{>0}:= \mathbb F_2\otimes_{\mathscr A} \mathcal P^{>0}_m,$ from which one has that $Q^{\otimes m} = (Q^{\otimes m})^0\,\bigoplus\, (Q^{\otimes m})^{>0}.$ For a polynomial $F\in \mathcal P_m,$ we denote by $[F]$ the classes in $Q^{\otimes m}$ represented by $F.$ If $\omega$ is a weight vector and $F\in \mathcal P_m^{\leq \omega},$ then denote by $[F]_\omega$ the classes in $(Q^{\otimes m})^{\omega}$ represented by $F.$ For a subset $\mathscr{C}\subset \mathcal P_m,$ we also write $|\mathscr C|$ for the cardinal of $\mathscr C$ and put $[\mathscr C] = \{[F]\, :\, F\in \mathscr C\}.$ If $\mathscr C\subset  P_m^{\leq \omega},$ then put $[\mathscr C]_{\omega} = \{[F]_{\omega}\, :\, F\in \mathscr C\}.$ Let us denote by $\mathscr {C}^{\otimes m}_n$ the set of all admissible monomials of degree $n$ in  $\mathcal P_m,$ 
%We put $$ \mathscr{C}^{0}_n := \mathscr {C}^{\otimes m}_n\cap (\mathcal P^{0}_m)_n,\ \mbox{and}\ \mathscr{C}^{>0}_n := \mathscr{C}^{\otimes m}_n\cap (\mathcal P^{>0}_m)_n.$$ 
and let $\omega$ be a weight vector of degree $n.$ By setting
 $$ \begin{array}{ll}
(\mathscr{C}^{\otimes m}_{n})^{\omega} := \mathscr {C}^{\otimes m}_n\cap \mathcal P_m^{\leq \omega},\ \ (\mathscr{C}^{\otimes m}_{n})^{\omega^{0}} := (\mathscr{C}^{\otimes m}_{n})^{\omega}\cap  \mathcal P^{0}_m,\ \ (\mathscr{C}^{\otimes m}_{n})^{\omega^{>0}} := (\mathscr{C}^{\otimes m}_{n})^{\omega}\cap  \mathcal P^{>0}_m,\\
(Q_n^{\otimes m})^{\omega^{0}}:= (Q^{\otimes m})^{\omega}\cap (Q^{\otimes m}_n)^{0},\ \ (Q_n^{\otimes m})^{\omega^{>0}} := (Q^{\otimes m})^{\omega}\cap (Q^{\otimes m}_n)^{>0},
\end{array}$$
then the sets $[(\mathscr{C}^{\otimes m}_{n})^{\omega}]_\omega,\, [(\mathscr{C}^{\otimes m}_{n})^{\omega^{0}}]_\omega$ and $[(\mathscr{C}^{\otimes m}_{n})^{\omega^{>0}}]_\omega$ are the bases of the $\mathbb F_2$-vector spaces $(Q_n^{\otimes m})^{\omega},\ (Q_n^{\otimes m})^{\omega^{0}}$ and $(Q_n^{\otimes m})^{\omega^{>0}},$ respectively.

\medskip

{\bf Main results and applications.} Let us now return to our study of the kernel of the Kameko homomorphism $(\widetilde {Sq^0_*})_{(5, 5(2^{1} - 1) + 18.2^{1})}$ and state our main results in greater detail. Firstly, by direct calculations using the results by Kameko \cite{M.K}, Singer \cite{W.S1}, Sum \cite{N.S1}, and  T\'in \cite{N.T}, we obtain the following, which is one of our main results and is crucial for an application on the dimension of $Q^{\otimes 6}.$

\begin{thm}\label{kq1} 
We have an isomorphism
$$ \mbox{\rm Ker}(\widetilde {Sq_*^0})_{(5, 5(2^{1} - 1) + 18.2^{1})} \cong  (Q_{5(2^{1} - 1) + 18.2^{1}}^{\otimes 5})^0  \bigoplus (Q_{5(2^{1} - 1) + 18.2^{1}}^{\otimes 5})^{\widetilde{\omega}^{>0}},$$
where $\widetilde{\omega} = (3,3,2,1,1)$ is the weight vector of the degree $5(2^{1} - 1) + 18.2^{1}.$
%If $X\in \mathscr {C}^{\otimes 5}_{5(2^{1} - 1) + 18.2^{1}}$ and $[X]\in \mbox{\rm Ker}(\widetilde {Sq_*^0})_{(5, 5(2^{1} - 1) + 18.2^{1})},$ then the weight vector $\widetilde{\omega}:=\omega(X)$ is one of the following sequences:
%$$ \widetilde{\omega}_{[1]} := (3,3,2,1,1), \ \  \widetilde{\omega}_{[2]} :=(3,3,2,3), \ \  \widetilde{\omega}_{[3]} :=(3,3,4,2).$$ 
\end{thm}

%Combining this and the above concepts, we see that 
%$$ \begin{array}{ll}
%\mbox{\rm Ker}(\widetilde {Sq_*^0})_{(5, 5(2^{1} - 1) + 18.2^{1})} &\cong (Q_{5(2^{1} - 1) + 18.2^{1}}^{\otimes 5})^0 \bigoplus \mbox{\rm Ker}(\widetilde {Sq_*^0})_{(5, 5(2^{1} - 1) + 18.2^{1})}\bigcap (Q_{5(2^{1} - 1) + 18.2^{1}}^{\otimes 5})^{>0}\\
%&\cong  (Q_{5(2^{1} - 1) + 18.2^{1}}^{\otimes 5})^0  \bigoplus\big(\bigoplus_{1\leq j\leq 3}(Q_{5(2^{1} - 1) + 18.2^{1}}^{\otimes 5})^{\widetilde{\omega}_{[j]}^{>0}}\big).
%\end{array}$$ 
\begin{rem}\label{nx2}
We are given in \cite{D.P6} that $(Q_{n}^{\otimes 5})^0  \cong \bigoplus_{1\leq s\leq 4}\bigoplus_{\ell(\mathcal J) = s}(Q_{n}^{\otimes\, \mathcal J})^{>0},$ where $$Q^{\otimes \mathcal J} = \langle [x_{j_1}^{t_1}x_{j_2}^{t_2}\ldots x_{j_s}^{t_s}]\;|\; t_i \in \mathbb N,\, i = 1,2, \ldots, s\}\rangle \subset Q^{\otimes 5}$$ with $\mathcal J= (j_1, j_2, \ldots, j_s),$  $1 \leq j_1 < \ldots < j_s \leq 5$,\, $1 \leq s \leq 4,$ and $\ell(\mathcal J):=s$ denotes the length of $\mathcal J.$ This implies that $ \dim((Q_{n}^{\otimes 5})^0)) = \sum_{1\leq  s\leq 4}\binom{5}{s}\dim((Q_{n}^{\otimes\, s})^{>0}), \ \mbox{for all $n\geq 0.$}$ On the other side, since $\xi(5(2^{1} - 1) + 18.2^{1}) = 3,$ by Peterson \cite{F.P} and Wood \cite{R.W}, the spaces $Q_{5(2^{1} - 1) + 18.2^{1}}^{\otimes 1}$ and $Q_{5(2^{1} - 1) + 18.2^{1}}^{\otimes 2}$ are trivial. Moreover, following  Kameko \cite{M.K} and Sum \cite{N.S1}, we have seen that $(Q_{5(2^{1} - 1) + 18.2^{1}}^{\otimes 3})^{>0}$ is $15$-dimensional and that $(Q_{5(2^{1} - 1) + 18.2^{1}}^{\otimes 4})^{>0}$ is $165$-dimensional. Therefore, we may conclude that $$\dim((Q^{\otimes 5}_{5(2^{1} - 1) + 18.2^{1}})^{0}= 15.\binom{5}{3} + 165.\binom{5}{4}=975.$$ 
%Moreover, $(Q_{5(2^{1} - 1) + 18.2^{1}}^{\otimes 5})^{\widetilde{\omega}_{[1]}^{0}}\cong (Q_{5(2^{1} - 1) + 18.2^{1}}^{\otimes 5})^0$ and $(Q_{5(2^{1} - 1) + 18.2^{1}}^{\otimes 5})^{\widetilde{\omega}_{[j]}^{0}} = 0,$ for $j = 2,\, 3.$ 
%Thus, to determine $\mbox{\rm Ker}(\widetilde {Sq_*^0})_{(5, 5(2^{1} - 1) + 18.2^{1})},$ we will compute $(Q_{5(2^{1} - 1) + 18.2^{1}}^{\otimes 5})^{\omega^{>0}}.$
\end{rem}

Next, due to Remarks \ref{nx1}, \ref{nx2}, and to Theorem \ref{kq1}, the space $Q^{\otimes 5}_{5(2^{1} - 1) + 18.2^{1}}$ will be determined by computing $(Q_{5(2^{1} - 1) + 18.2^{1}}^{\otimes 5})^{\omega^{>0}}.$ To accomplish this, we use the method described above to explicitly indicate all the admissible monomials in the set $(\mathscr{C}^{\otimes 5}_{5(2^{1}-1)  +18.2^{1}})^{\widetilde{\omega}^{>0}}.$ As a result, it reads as follows.

\begin{thm}\label{kq2}
There exist exactly $925$ admissible monomials of degree $5(2^{1} - 1) + 18.2^{1}$ in $\mathcal P_5^{>0}$ such that their weight vectors are $\widetilde{\omega}.$ Consequently, $(Q_{5(2^{1} - 1) + 18.2^{1}}^{\otimes 5})^{\widetilde{\omega}^{>0}}$ has dimension $925.$

%\begin{thm}\label{kq2}
%The following statements are true:
%\begin{enumerate}
%\item[{\rm i)}] The $\mathbb F_2$-vector spaces $(Q_{5(2^{1} - 1) + 18.2^{1}}^{\otimes 5})^{\widetilde{\omega}_{[j]}^{>0}}$ are trivial, for $j = 2,\, 3.$

%\item[{\rm ii)}] $(Q_{5(2^{1} - 1) + 18.2^{1}}^{\otimes 5})^{\widetilde{\omega}_{[1]}^{>0}}$ is $925$-dimensional.
%\end{enumerate}

\end{thm}

This theorem, together with the fact that $Q_{5(2^{t}-1)+18.2^{t}}^{\otimes 5} = (Q_{5(2^{t}-1)+18.2^{t}}^{\otimes 5})^0\,\bigoplus\, (Q_{5(2^{t}-1)+18.2^{t}}^{\otimes 5})^{>0},$ yields an immediate corollary that

\begin{corl}\label{hq}
The space $Q_{5(2^{t}-1)+18.2^{t}}^{\otimes 5}$ is $730$-dimensional if $t = 0,$ and is $2630$-dimensional if $t\geq 1.$
\end{corl}

\medskip

%Very recently, an useful hypothesis to study hit problems is given by Sum \cite{N.S2}, who sets up a relation between the minimal set of $\mathscr A$-generators for the polynomial algebras $\mathcal P_m$ and $\mathcal P_{m-1}$ in each the certain degree. We refer the reader to section three for a more detailed treatment. Calculations by Peterson \cite{F.P}, Kameko \cite{M.K} and Sum \cite{N.S1} showed that this prediction holds for $m\leq 4.$  By Sum himself  \cite{N.S2} and our previous results \cite{D.P1, D.P2, D.P3}, it is known to hold in the case $m = 5$ and some certain degrees. Due to the proof of Theorem \ref{kq2} and Remark \textcolor[rgb]{0.0,0.0,1.0}{4.2.4} in section four, the conjecture is also true for $m=5$ and the generic degree $5(2^{t} - 1) + 18.2^{t}.$ 

%\medskip

As applications, one would also be interested in applying results and techniques of hit problems into the cases of higher ranks $m$ of $Q^{\otimes m}$ and the modular representations of the general linear groups (see also the relevant discussions in  literatures  \cite{J.B}, \cite{N.M1, N.M2}, \cite{T.N2}, \cite{W.W, W.W2}). Two applications below of the contributions of this paper are also not beyond this target.
\medskip

{\bf First application: the dimension of \boldmath{$Q^{\otimes 6}.$}} The hit problem of six variables has been not yet known. Using Corollary \ref{hq} for the case $t \geq 1$ and a result in Sum \cite{N.S1}, we state that

\begin{thm}\label{kqbs}
With the generic degree $5(2^{t+4} - 1) + 41.2^{t+4},$ where $t$ an arbitrary positive integer, then the $\mathbb F_2$-vector space $Q^{\otimes 6}$ has dimension $165690$ in this degree.
\end{thm}

Observing from Corollary \ref{hq} and Theorem \ref{kqbs}, the readers can notice that the dimensions of $Q^{\otimes 5}$ and $Q^{\otimes 6}$ in degrees given are very large. So, a general approach to hit problems, other than providing a mononial basis of the vector space $Q_n^{\otimes m},$  is to find upper/lower bounds on the dimension of this space. However, in this work, we have not studied this side of the problem and it is our concern the next time. It is remarkable that, we have Kameko's conjecture \cite{M.K} on an upper bound for the dimension of $Q_n^{\otimes m},$ but unfortunately, it was refuted for $m\geq 5$ by the brilliant work of Sum \cite{N.S}.

\medskip

{\bf Second application: the behavior of the fifth Singer transfer}. We adopt Corollary \ref{hq} for $t = 0,$ together with a fact of the Adams $E^{2}$-term, ${\rm Ext}_{\mathscr A}^{5, 5+*}(\mathbb F_2, \mathbb F_2)$, to obtain information about the behavior of  Singer's cohomological transfer in the bidegree $(5, 5+(5(2^{0}-1) + 18.2^{0}))$. More precisely, it is known, the calculations of Lin \cite{W.L}, and Chen \cite{T.C} imply that ${\rm Ext}_{\mathscr A}^{5, 5+(5(2^{t}-1)  +18.2^{t})}(\mathbb F_2, \mathbb F_2) = \langle h_tf_t \rangle$ and  $h_tf_t = h_{t+1}e_t \neq 0$ for all $t\geq 0.$ So, to determine the transfer map in the above bidegree, we shall compute the dimension of (the domain of the fifth transfer) $(\mathbb F_2 \otimes_{GL_5}P_{\mathscr A}((\mathcal P_5)^{*}))_{5(2^{0}-1)  +18.2^{0}}$ by using a mononial basis of $Q^{\otimes 5}_{5(2^{0}-1) + 18.2^{0}}.$ (We emphasize that computing the domain of $Tr_m^{\mathscr A}$ in each degree $n$ is very difficult, particularly for values of $m$ as large as $m = 5.$ The understanding of special cases should be a helpful step toward the solution of the general problem. Moreover, we believe, in principle, that our method could lead to a full analysis of $\mathbb F_2 \otimes_{GL_m}P_{\mathscr A}((\mathcal P_m)^{*})$ in each $m$ and degree $n > 0$, as long as nice decompositions of the space of $GL_m$-invariants of $Q^{\otimes m}$ in degrees given. However, the difficulty of such a task must be monumental, as $Q^{\otimes m}$ becomes much larger and harder to understand with increasing $m.$) Details for this application are as follows. It may need to be recalled that by the previous discussions \cite{D.P3}, we get the technical proposition below.

\begin{prop}\label{md18}
The following hold:
\begin{enumerate}
\item[{\rm i)}] If $Y\in \mathscr C_{5(2^{0}-1)  +18.2^{0}}^{\otimes 5},$ then $\overline{\omega}:=\omega(Y)$ is one of the following sequences:
$$ \overline{\omega}_{[1]} := (2,2,1,1), \ \   \overline{\omega}_{[2]} :=(2,2,3), \ \   \overline{\omega}_{[3]} :=(2,4,2),$$ 
$$  \overline{\omega}_{[4]} :=(4,1,1,1), \ \   \overline{\omega}_{[5]} :=(4,1,3), \ \   \overline{\omega}_{[6]} :=(4,3,2).$$

\item[{\rm ii)}] 
$|(\mathscr{C}^{\otimes 5}_{5(2^{0}-1)  +18.2^{0}})^{\overline{\omega}_{[k]}}|  = \left\{\begin{array}{ll}
300 &\mbox{{\rm if}}\ k= 1,\\
15 &\mbox{{\rm if}}\ k= 2, 5,\\
10 &\mbox{{\rm if}}\ k= 3,\\
110 &\mbox{{\rm if}}\ k= 4,\\
280 &\mbox{{\rm if}}\ k= 6.
\end{array}\right.$
\end{enumerate}
\end{prop}
One should note that $|(\mathscr{C}^{\otimes 5}_{5(2^{0}-1)  +18.2^{0}})^{\overline{\omega}_{[k]}}| = |(\mathscr{C}^{\otimes 5}_{5(2^{0}-1)  +18.2^{0}})^{\overline{\omega}_{[k]}^{>0}}|$ for $k = 2, 3$, and that $|(\mathscr{C}^{\otimes 5}_{5(2^{0}-1)  +18.2^{0}})^{\overline{\omega}_{[2]}^{0}}| = 0 = |(\mathscr{C}^{\otimes 5}_{5(2^{0}-1)  +18.2^{0}})^{\overline{\omega}_{[3]}^{0}}|.$ Moreover, $\dim(Q^{\otimes 5}_{5(2^{0}-1)  +18.2^{0}}) = \sum_{1\leq k\leq 6}|(\mathscr{C}^{\otimes 5}_{5(2^{0}-1)  +18.2^{0}})^{\overline{\omega}_{[k]}}| = 730.$ Next, applying these results, we explicitly compute the subspaces of $GL_5$-invariants $((Q_{5(2^{0} - 1) + 18.2^{0}}^{\otimes 5})^{\overline{\omega}_{[k]}})^{GL_5},$ for $1\leq k\leq 6,$ and obtain

\begin{thm}\label{kq3} 
The following assertions are true:
\begin{enumerate}
\item[{\rm i)}] $((Q_{5(2^{0} - 1) + 18.2^{0}}^{\otimes 5})^{\overline{\omega}_{[k]}})^{GL_5} = 0$ with  $k \in \{1, 2, 3, 5, 6\}.$

\item[{\rm ii)}] $((Q_{5(2^{0} - 1) + 18.2^{0}}^{\otimes 5})^{\overline{\omega}_{[4]}})^{GL_5} = \langle [\Re'_4]_{\overline{\omega}_{[4]}} \rangle,$ where
$$ \begin{array}{ll}
\medskip
\Re'_4&= x_1x_2x_3x_4x_5^{14} + x_1x_2x_3x_4^{14}x_5 + x_1x_2x_3^{14}x_4x_5 + x_1x_2^{3}x_3x_4x_5^{12}\\
\medskip
&\quad + x_1x_2^{3}x_3x_4^{12}x_5 + x_1x_2^{3}x_3^{12}x_4x_5 + x_1^{3}x_2x_3x_4x_5^{12} + x_1^{3}x_2x_3x_4^{12}x_5\\
\medskip
&\quad + x_1^{3}x_2x_3^{12}x_4x_5 + x_1^{3}x_2^{5}x_3x_4x_5^{8} + x_1^{3}x_2^{5}x_3x_4^{8}x_5 + x_1^{3}x_2^{5}x_3^{8}x_4x_5. 
\end{array}$$
\end{enumerate}
\end{thm}

Now, because $(\mathbb F_2 \otimes_{GL_5}P_{\mathscr A}((\mathcal P_5)^{*}))_{5.(2^{0}-1) + 18.2^{0}}$ is isomorphic to $(Q^{\otimes 5}_{5.(2^{0}-1) + 18.2^{0}})^{GL_5},$ by Theorem \ref{kq3}, we have the following estimate:
$$ \begin{array}{ll}
 \dim (\mathbb F_2 \otimes_{GL_5}P_{\mathscr A}((\mathcal P_5)^{*}))_{5.(2^{0}-1) + 18.2^{0}} &= \dim (Q^{\otimes 5}_{5.(2^{0}-1) + 18.2^{0}})^{GL_5}\nonumber\\
& \leq \sum_{1\leq k\leq 6}\dim((Q_{5(2^{0} - 1) + 18.2^{0}}^{\otimes 5})^{\overline{\omega}_{[k]}})^{GL_5}\leq 1.
\end{array}$$
On the other side, as shown in section one, $\{h_t|\, t\geq 0\}\subset {\rm Im}(Tr^{\mathscr A}_1),$ and $\{f_{t}|\, t\geq 0\}\subset {\rm Im}(Tr^{\mathscr A}_4).$ Combining this with the fact that the total transfer $\bigoplus_{m\geq 0} Tr_m^{\mathscr A}$ is a homomorphism of algebras, it may be concluded that the non-zero element $h_tf_t \in {\rm Ext}_{\mathscr {A}}^{5, 23.2^{t}}(\mathbb F_2,\mathbb F_2) $ is in the image of $Tr^{\mathscr A}_5$ for all $t\geq 0.$ This could be directly proved as in Appendix. This statement implies that $$ \dim (\mathbb F_2 \otimes_{GL_5}P_{\mathscr A}((\mathcal P_5)^{*}))_{5.(2^{0}-1) + 18.2^{0}}\geq 1,$$
and therefore $(\mathbb F_2 \otimes_{GL_5}P_{\mathscr A}((\mathcal P_5)^{*}))_{5.(2^{0}-1) + 18.2^{0}}$ is one-dimensional. As a consequence, we immediately obtain

 \begin{corl}\label{hq2}
The cohomological transfer
$$Tr^{\mathscr A}_5: (\mathbb F_2 \otimes_{GL_5}P_{\mathscr A}((\mathcal P_5)^{*}))_{5(2^0-1) + 18.2^0}\to {\rm Ext}_{\mathscr {A}}^{5, 5+5(2^0-1) + 18.2^0}(\mathbb F_2,\mathbb F_2)$$
is an isomorphism. Consequently, Conjecture \ref{gtSin} holds in the rank 5 case and the degree $5(2^{0}-1) + 18.2^{0}.$ 
\end{corl}

\medskip

{\bf Comments and open issues.} From the above results, it would be interesting to see that $Q^{\otimes 5}$ is $730$-dimensional in degree $5(2^{0} - 1) + 18.2^{0},$ but the space of $GL_5$-coinvariants of it in this degree is only one-dimensional. In general, it is quite efficient in using the results of the hit problem of five variables to study $\mathbb F_2 \otimes_{GL_5}P_{\mathscr A}((\mathcal P_5)^{*}).$ This provides a valuable method for verifying Singer's open conjecture on the fifth algebraic transfer. We now close the introduction by discussing about Conjecture \ref{gtSin} in the rank 5 case and the internal degree $n_t:=5(2^{t}-1) + 18.2^{t}$ for all $t \geq 1.$ Let us note again that the iterated Kameko homomorphism $((\widetilde {Sq^0_*})_{(5, n_t)})^{t-1}:  Q^{\otimes 5}_{n_t} \to Q^{\otimes 5}_{n_1}$ is an $\mathbb F_2GL_5$-module isomorphism for all $t\geq 1.$ So, from a fact of ${\rm Ext}_{\mathscr {A}}^{5, 5 + n_1}(\mathbb F_2,\mathbb F_2)$, to check Singer's conjecture  in the above degree, we need only determine $GL_5$-coinvariants of $Q^{\otimes 5}_{n_t}$ for $t = 1.$ We must recall that Kameko's map $(\widetilde {Sq^0_*})_{(5, n_1)}:  Q^{\otimes 5}_{n_1} \to Q^{\otimes 5}_{n_0}$ is an epimorphism of $GL_5$-modules. On the other side, as shown before, the non-zero element $h_1f_1\in {\rm Ext}_{\mathscr {A}}^{5, 5+n_1}(\mathbb F_2,\mathbb F_2)$ is detected by the fifth transfer. From these data and Theorem \ref{kq3}, one has an estimate $$0\leq \dim((\mathbb F_2 \otimes_{GL_5}P_{\mathscr A}((\mathcal P_5)^{*}))_{n_1}) -1 \leq \dim ({\rm Ker}(\widetilde {Sq^0_*})_{(5,n_1)})^{GL_5}.$$ Moreover, basing the proof of Theorem \ref{kq3} together with a few simple arguments, it follows that the elements in $(\mathbb F_2 \otimes_{GL_5}P_{\mathscr A}((\mathcal P_5)^{*}))_{n_1}$ are dual to the classes
$$ \begin{array}{ll}
\medskip
&\gamma[x_1^{3}x_2^{3}x_3^{3}x_4^{3}x_5^{29} + x_1^{3}x_2^{3}x_3^{3}x_4^{29}x_5^{3} + x_1^{3}x_2^{3}x_3^{29}x_4^{3}x_5^{3} + x_1^{3}x_2^{7}x_3^{3}x_4^{3}x_5^{25}+ x_1^{3}x_2^{7}x_3^{3}x_4^{25}x_5^{3}  \\
\medskip
&\quad + x_1^{3}x_2^{7}x_3^{25}x_4^{3}x_5^{3} + x_1^{7}x_2^{3}x_3^{3}x_4^{3}x_5^{25} + x_1^{7}x_2^{3}x_3^{3}x_4^{25}x_5^{3}+ x_1^{7}x_2^{3}x_3^{25}x_4^{3}x_5^{3} + x_1^{7}x_2^{11}x_3^{3}x_4^{3}x_5^{17} \\
&\quad + x_1^{7}x_2^{11}x_3^{3}x_4^{17}x_5^{3} + x_1^{7}x_2^{11}x_3^{17}x_4^{3}x_5^{3}] + [\zeta],
\end{array}$$
where $\gamma\in \mathbb F_2,$ and $[\zeta]\in \text{Ker}(\widetilde {Sq^0_*})_{(5, n_1)}.$ It could be noticed that calculating explicitly these elements is not easy.  However, in view of our previous works \cite{D.P2, D.P6}, and motivated by the above computations, we have the following prediction.

\begin{gt}\label{gtP}
For each $t\geq 1,$ the space of $GL_5$-invariants elements of ${\rm Ker}(\widetilde {Sq^0_*})_{(5,n_t)}$ is trivial.  Consequently, the coinvariant $(\mathbb F_2 \otimes_{GL_5}P_{\mathscr A}((\mathcal P_5)^{*}))_{n_t}$ is 1-dimensional.
\end{gt}

Since $h_{t}f_{t}\in {\rm Im}(Tr_5^{\mathscr A}),$ for all $t\geq 0,$ if Conjecture \ref{gtP} is true, then $Tr_5^{\mathscr A}$ is also isomorphism when acting on the coinvariant $(\mathbb F_2 \otimes_{GL_5}P_{\mathscr A}((\mathcal P_5)^{*}))_{n_t}$ for  $t\geq 1,$ and so, Conjecture \ref{gtSin} holds in bidegree $(5, 5+n_t)$. We also wish that our predictions are correct. If not, Singer's conjecture will be disproved.  We leave these issues as future research. At the same time, we also appreciate that some readers may have an interest in solving them.

\medskip

{\bf Overview.} Let us give a brief outline of the contents of this paper. Section three contains a brief review of Steenrod squares and some useful linear transformations. The dimensions of the polynomial algebras $\mathcal P_5$ and $\mathcal P_6$ in the generic degrees $n_t = 5(2^{t}-1) + 18.2^{t}$ and $5(2^{t+4} - 1) + n_1.2^{t+4}$ are respectively obtained in section four by proving Theorems \ref{kq1}, \ref{kq2}, and \ref{kqbs}. Section five is to present the proof of Theorem \ref{kq3}. In the remainder of the text, we give a direct proof of an event claimed above that the non-zero elements $h_tf_t \in {\rm Ext}_{\mathscr {A}}^{5, 23.2^{t}}(\mathbb F_2,\mathbb F_2) $ are detected by $Tr_5^{\mathscr A}$. The proof is based on a representation in the lambda algebra of the fifth Singer transfer. Finally, we describe the set  $(\mathscr C_{n_1}^{\otimes 5})^{\widetilde{\omega}^{>0}}$ and list some the admissible monomials in $\mathscr C_{n_0}^{\otimes 5}$ and the strictly inadmissible monomials in ($\mathcal P_5^{>0})_{n_1}.$ 
%Finally, we provide an algorithm in MAGMA for verifying the dimension of $Q^{\otimes 5}$ in degree $n_t$ and studying the hit problem in general.

\begin{acknow}
The author is indebted to the anonymous referees for very useful comments and suggestions on the preceding versions of this work.
%The author would like to give my deepest sincere thanks to Professor N. Sum (Quy Nhon University, Sai Gon University), my thesis advisor, for meaningful discussions. 
%First of all, I gratefully acknowledge the anonymous referee for reading the paper carefully and providing thoughtful comments, many of which have resulted in changes to the revised version of the manuscript.
%The author would like to express his deep gratitude to Prof. N. Sum, his thesis advisor, for having proposed this research topic,  and having devoted a great deal of time to discuss details of the research with the author.  
I gratefully acknowledge Professor W. Singer for many enlightening e-mail exchanges. 
\end{acknow}

\section{The Necessary Preliminaries}\label{s2}
This section begins with a few words on the Steenrod algebra over $\mathbb F_2$ and ends with a brief sketch of some homomorphisms in \cite{N.S1}. At the same time, we prove some elementary results that will be used in the rest of this text.

%In this section, we provide a brief overview of the Steenrod algebra over $\mathbb F_2$ and its application. At the same time, we review some concepts and the relevant results in \cite{P.S1, N.S1}, which will be used in the next section. 
%For the concepts and results related to the Peterson hit problem are not shown in this paper, we will refer to \cite{D.P3, N.S1} for details. 

\subsection{Steenrod squares and their properties}\label{sub1-1}
The mod 2 Steenrod algebra $\mathscr A$ was defined by Cartan \cite{Cartan} to be the algebra of stable cohomology operations for mod 2 cohomology. This algebra is generated by the Steenrod squares $Sq^i: H^{n}(X, \mathbb{F}_2) \to H^{n+i}(X, \mathbb{F}_2),$ for $i\geq 0,$ where $H^{n}(X, \mathbb{F}_2)$ denotes the $n$-th singular cohomology group of a topological space $X$ with coefficient over $\mathbb F_2.$ Steenrod and Epstein \cite{S.E} showed that these squares are characterized by the following 5 axioms:

\begin{enumerate}
\item[(i)] $Sq^i$ is an additive homomorphism and is natural with respect to any $f: X\to Y.$ So $f^{*}(Sq^{i}(x)) = Sq^{i}(f^{*}(x)).$ 
\item[(ii)] $Sq^0$ is the identity homomorphism.
\item[(iii)] $Sq^{i}(x) = x\smile x$ for all $x\in H^{i}(X, \mathbb F_2)$ where $\smile$ denotes the \textit{cup product} in the graded-commutative ring $H^*(X, \mathbb{F}_2).$ 
\item[(iv)]  If $i > \deg(x),$ then $Sq^{i}(x) = 0.$
\item[(v)] Cartan's formula: $Sq^n(x\smile y) =  \sum_{i+j = n}Sq^i(x)\smile Sq^{j}(y).$

In addition, Steenrod squares have the following properties:

$\bullet$ $Sq^1$ is the Bockstein homomorphism of the coefficient sequence: $0\to \mathbb Z/2\to \mathbb Z/4\to \mathbb Z/2\to 0.$

$\bullet$ $Sq^{i}$ commutes with the connecting morphism of the long exact sequence in cohomology. In particular, it commutes with respect to suspension $H^{n}(X, \mathbb F_2)\cong H^{n+1}(\Sigma X, \mathbb F_2).$

$\bullet$ They satisfy the Adem relations:
$Sq^iSq^j  =  \sum_{0\leq t\leq [i/2]}\binom{j-t-1}{i-2t}Sq^{i+j-t}Sq^t,\,\,0 < i < 2j,$ where the binomial coefficients are to be interpreted mod 2. These relations, which were conjectured by Wu \cite{Wu} and established by Adem \cite{Adem}, allow one to write an arbitrary composition of Steenrod squares as a sum of Serre-Cartan basis elements.
\end{enumerate}

Note that the structure of the cohomology $H^*(X, \mathbb{F}_2)$ is not only as graded commutative $\mathbb F_2$-algebra, but also as an $\mathscr A$-module. In many cases, the $\mathscr A$-module structure on $H^{*}(X, \mathbb F_2)$ provides additional information on $X.$ 
%For instance, one has the following.

\end{document}